\documentclass[12pt]{article}
\usepackage{latexsym}
\usepackage{rotating}
\usepackage{amssymb}
\usepackage{amsmath}
\usepackage{amsthm}
\usepackage{enumerate}
\usepackage{graphicx}
\usepackage{comment}
\usepackage{epsfig}

\usepackage{amsmath,amsfonts,amscd}

\newtheorem{theorem}{Theorem}[section]
\newtheorem{proposition}{Proposition}[section]
\newtheorem{lemma}{Lemma}[section]
\newtheorem{corollary}{Corollary}[section]

\newcommand{\psiAm}[2]{\psi_{#1,#2}^{A,-}}
\newcommand{\Pbr}{P_{\psi_r}}
\newcommand{\Abr}{A_{\psi_r}}
\newcommand{\Bbr}{B_{\psi_r}}
\newcommand{\Px}{P_{\chi}}
\newcommand{\Ax}{A_{\chi}}
\newcommand{\Bx}{B_{\chi}}

\theoremstyle{definition}
\newtheorem{definition}{Definition}[section]
\newtheorem{remark}{Remark}[section]
\newtheorem{example}{Example}[section]
\newcommand{\BbZ} {\mathbb Z}
\newcommand{\BbN} {\mathbb N}
\newcommand{\Q} {\mathbb Q}
\newcommand{\U} {\mathcal{U}}
\newcommand{\Ut} {\tilde{\mathcal{U}}}

\newcommand{\ddb}{d_\beta(1)}
\newcommand{\ddbs}{d^*_{\beta}(1)}

\title{On univoque Pisot numbers}

\author{Jean-Paul Allouche\thanks{Research partially supported by MENESR,
ACI NIM 154 Num\'eration.} \\
CNRS, LRI, B\^atiment 490 \\
Universit\'e Paris-Sud \\
91405 Orsay Cedex, France \\
{\tt allouche@lri.fr}
\and
Christiane Frougny \\
LIAFA, CNRS UMR 7089 \\
2 place Jussieu \\
75251 Paris Cedex 05, France \\
and Universit\'{e} Paris 8 \\
{\tt Christiane.Frougny@liafa.jussieu.fr}
\and
\setcounter{footnote}{6}%
Kevin G. Hare\thanks{Research supported, in part by NSERC of Canada.} \\
Department of Pure Mathematics \\
University of Waterloo \\
Waterloo, Ontario, Canada, N2L 3G1 \\
{\tt kghare@math.uwaterloo.ca}
}

\begin{document}

\maketitle

\begin{abstract}
We study Pisot numbers $\beta \in (1, 2)$ which are univoque, i.e., such that 
there exists only one representation of $1$ as 
    $1 = \sum_{n \geq 1} s_n\beta^{-n}$, with
$s_n \in \{0, 1\}$. We prove in particular that there exists a smallest univoque
Pisot number, which has degree $14$. Furthermore we give the smallest limit point
of the set of univoque Pisot numbers. 
\end{abstract}



\section{Introduction}

Representations of real numbers in non-integer bases were introduced 
by R\'enyi \cite{Re} and first studied by R\'enyi and by Parry \cite{Pa, Re}.
Among the questions that were addressed is the  uniqueness of 
representations. 
Given a sequence $(s_n)_{n\geq 1}$, 
Erd\H{o}s, Jo\'o and Komornik, \cite{EJK}, gave  
a purely combinatorial characterization for when there exists $\beta \in (1,2)$
such that $1 = \sum_{n \geq 1} s_n\beta^{-n}$ is the unique 
representation of 1.  This set of
binary sequences is essentially the same as a set studied by Cosnard and the 
first author \cite{All, AC1, AC3} in the context of iterations of unimodal 
continuous maps of the unit interval.

Following \cite{KL, KLP}, a number $\beta>1$ is said to be {\em univoque} 
if there exists a unique sequence of integers $(s_n)_{n \ge 1}$, 
with $0 \le s_n <\beta$, such that $1=\sum_{n \ge 1}s_n \beta^{-n}$. 
(Note that we consider only the representation of $1$. The uniqueness of the 
representation of real numbers in general was studied in particular in~\cite{GS}.)
Using the characterization of \cite{EJK}, Komornik and Loreti constructed in 
\cite{KL} the smallest real number in $(1, 2)$ for which $1$ has a unique 
representation. Its representation happens to be the famous Thue-Morse sequence 
(see for example \cite{AS}).

Are there univoque Pisot numbers? It is worth noting that if the base $\beta$ is the 
``simplest'' non-integer Pisot number, i.e., the golden ratio, then the number $1$ 
has infinitely many representations. In this paper we study the univoque Pisot numbers 
belonging to $(1, 2)$. We prove in particular (Theorem~\ref{thm:finite2}) that there 
exists a smallest univoque Pisot number, and we give explicitly the least three 
univoque Pisot numbers in $(1,2)$: they are the roots in $(1, 2)$ of the polynomials
$$
\begin{array}{lll}
& {x}^{14}-2{x}^{13}+{x}^{11}-{x}^{10}-{x}^{7}+{x}^{6} -{x}^{4}+{x}^{3}-x+1
 & (\mathrm{root} \approx 1.8800), \\
& {x}^{12}-2{x}^{11}+{x}^{10}-2{x}^{9}+{x}^{8}-{x}^{3}+{x}^{2}-x+1, &
   (\mathrm{root} \approx 1.8868), \\
& x^4 - x^3 - 2x^2 + 1, &
    (\mathrm{root} \approx 1.9052).
\end{array}
$$
The last number is the smallest limit point of the set of univoque Pisot numbers
(Theorem~\ref{thm:chi}). We also prove that $2$ is a limit point of univoque Pisot numbers.

\section{Definitions and reminders}
\subsection{Infinite words}
Let ${\mathbb N}_+$ denote the set of positive integers.
Let $A$ be a finite alphabet. We define $A^{{\mathbb N}_+}$ to be the set 
of infinite sequences (or infinite words) on $A$:
$$
A^{{\mathbb N}_+} := \{s=(s_n)_{n \ge 1} \ | \ \forall n \ge 1, \ s_n \in A\}.
$$
This set is equipped with the distance $\rho$ defined by: 
if $s=(s_n)_{n \ge 1}$ and $v=(v_n)_{n \ge 1}$ belong to $A^{{\mathbb N}_+}$, 
then $\rho(s,v) := 2^{-r}$ if $s \neq v$ and $r := \min\{n \mid s_n \neq v_n\}$, 
and $\rho(s,v)=0$ if $s=v$.
The topology on the set $A^{{\mathbb N}_+}$ is then the product topology, 
and it makes $A^{{\mathbb N}_+}$ a compact metric space.

A sequence $(s_n)_{n \geq 1}$ in $A^{{\mathbb N}_+}$ is said to be
{\it periodic\,} if there exists an integer $T \geq 1$, called a {\em period}
of the sequence, such that $s_{n+T} = s_n$ for all $n \geq 1$. A sequence
$(s_n)_{n \geq 1}$ in $A^{{\mathbb N}_+}$ is said to be {\it eventually
periodic\,} if there exists an integer $n_0 \geq 0$ such that the sequence
$(s_{n+n_0})_{n \geq 1}$ is periodic.

If $w$ is a (finite) word, we denote by $w^{\infty}$ the infinite word obtained
by concatenating infinitely many copies of $w$ (this is in particular a periodic
sequence, and the length of $w$, usually denoted by $|w|$, is a period).

\subsection{Base $\beta$ representations}

Let $\beta$ be a real number $>1$. A {\em $\beta$-representation}
of the real number $x \in [0,1]$ is an infinite sequence of integers
$(x_n)_{n \ge 1}$ such that $x = \sum_{n \ge 1}x_n \beta^{-n}$. 
If a representation ends in infinitely many zeros, say, is of the form $w0^\infty$,
then the ending zeros are omitted and the representation is said to 
be {\em finite}.
The reader is referred to \cite[Chapter 7]{Lo} for more on these topics.

\subsubsection{Greedy representations}
\label{sec:greedy}

A special representation of a number $x$, called the {\em greedy 
$\beta$-expansion}, is the infinite sequence $(x_n)_{n \ge 1}$ obtained by 
using the greedy algorithm of R\'enyi~\cite{Re}.

Denote by $\lfloor y \rfloor$ and $\{y\}$ the integer part 
and the fractional part of the real number $y$. Set $r_0 := x$ and, for $n \ge 1$, 
let $x_n := \lfloor \beta r_{n-1} \rfloor$, $r_n := \{\beta r_{n-1}\}$.
Then $x= \sum_{n \ge 1} x_n \beta^{-n}$.

Intuitively, the digit $x_n$ is chosen so that it is the maximal choice allowed
    at each step. The digits $x_n$ obtained by the greedy 
algorithm belong to the alphabet $A_\beta = \{0,1,\ldots,\lfloor \beta \rfloor\}$
if $\beta$ is not an integer, which will always be the case in this work.
It is clear from the definition that amongst the $\beta$-representations of 
a number, the greedy $\beta$-expansion is the largest in the lexicographic 
order (denoted by $\leq_{lex}$ and $<_{lex}$).
The greedy $\beta$-expansion of $x$ will be denoted by
$d_\beta(x) := (x_n)_{n \ge 1}$. 

The greedy $\beta$-expansion of $1$ plays an important role. Set $\ddb =(e_n)_{n \ge 1}$ 
and define 
$$
\ddbs := 
\left\{
\begin{array}{ll}
\ddb \ &\mbox{\rm if $\ddb$ is infinite} \\
(e_1 \cdots e_{m-1} (e_m - 1))^\infty \ 
&\mbox{\rm if $\ddb=e_1 \cdots e_{m-1}e_m$ is finite}.
\end{array}
\right.
$$
Of course if $\ddb$ is finite, the sequence $\ddbs$ is also a $\beta$-representation 
of $1$.

\noindent
Denote by $\sigma$ the shift on $A_\beta^{{\mathbb N}_+}$: for any sequence
$s = (s_n)_{n \geq 1}$ in $A_\beta^{{\mathbb N}_+}$, the sequence 
$v = \sigma(s)$ is defined by  $v = (v_n)_{n \geq 1} := (s_{n+1})_{n \geq 1}$.
We recall some useful results.

\begin{theorem}\label{Pa} {\rm \cite{Pa}} 
Let $s = (s_n)_{n \geq 1}$ be a sequence in $A_{\beta}^{{\mathbb N}_+}$. Then

\begin{itemize}

\item the sequence $s$ is the greedy $\beta$-expansion of
some $x \in [0,1)$ 
if and only if
$$
\forall k \geq 0, \ \ \sigma^k(s) <_{lex} \ddbs
$$

\item the sequence $s$ is the greedy $\beta$-expansion of\; $1$
 for some $\beta >1$ if and only if
$$
\forall k \geq 1, \ \ \sigma^k(s) <_{lex} s.
$$

\end{itemize}
\end{theorem}

\subsubsection{Lazy representations}
\label{sec:lazy}

Another distinguished $\beta$-representation of the real number $x$ 
is the 
so-called {\em lazy} expansion, which is the smallest in the lexicographic 
order among the $\beta$-representations of $x$ on the alphabet $A_\beta$.
Denote by $\ell_\beta(x)=(x_n)_{n \geq 1} $ the lazy $\beta$-expansion of $x$.

To compute it, intuitively we have to choose $x_n$ to be as small as possible at each
    step. The algorithm to obtain the lazy expansion is the following.
    Let $B := \sum_{n \geq 1} \frac{\lfloor \beta \rfloor}{\beta^n}=\frac{\lfloor \beta \rfloor}{\beta -1}$.
Set $r_0 := x$ and, 
    for $n \ge 1$, let $x_n := \max(0,\lceil \beta r_{n-1}-B \rceil)$, 
    $r_n := \beta r_{n-1} - x_n$.
Then $x= \sum_{n \ge 1} x_n \beta^{-n}$, where the $(x_n)$ forms the lazy
    $\beta$-expansion.

Let $s=(s_n)_{n \ge 1}$ be in $A_\beta^{{\mathbb N}_+}$. Denote by 
$\overline{s_n} := \lfloor \beta  \rfloor-s_n$ the ``complement'' of $s_n$,
and by extension $\bar{s} := (\overline{s_n})_{n \ge 1}$.
Then the following characterization of lazy expansions holds true.

\begin{theorem}\label{EJK1} {\rm \cite{EJK,DK}} 
Let $s = (s_n)_{n \geq 1}$ be a sequence in $A_{\beta}^{{\mathbb N}_+}$. Then

\begin{itemize}

\item the sequence $s$ is the lazy $\beta$-expansion of some $x \in [0,1)$ 
if and only if
$$
\forall k \ge 0, \ \ \sigma^k(\bar{s}) <_{lex} \ddbs
$$
\item the sequence $s$ is the lazy $\beta$-expansion of\, $1$ 
for some $\beta >1$ if and only if
$$
\forall k \ge 1, \ \ \sigma^k(\bar{s}) <_{lex} s.
$$
\end{itemize}
\end{theorem}

\begin{example}\label{GR}
Take $\psi_1=\frac{1+\sqrt{5}}{2}$ the golden ratio.
The greedy $\beta$-expansion of\, $1$ is $d_{\psi_1}(1)=11$,
$d_{\psi_1}^*(1)=(10)^\infty$, and the lazy expansion of\, $1$ 
is $\ell_{\psi_1}(1)=01^\infty$.
\end{example}

\subsection{Univoque real numbers}

Following \cite{KL, KLP}, a number $\beta > 1$ is said 
to be {\em univoque} if there exists a unique sequence of integers $(s_n)_{n \ge 1}$, 
with $0 \le s_n <\beta$, such that $1=\sum_{n \ge 1}s_n \beta^{-n}$.
In this case the sequence $(s_n)_{n \ge 1}$ coincides both with the greedy 
and with the lazy $\beta$-expansion of $1$.
Remark that the number 2 is univoque, but we will be concerned
with non-integer real numbers in this paper.

Note that some authors call ``univoque'' the real numbers $x$ having
a unique $\beta$-representation (see \cite{DK2}). Binary sequences 
$(s_n)_{n \geq 1}$ such that the convergent sum 
$\sum_{n \geq 1}{s_n}\beta^{-n}$ uniquely determines the sequence 
$(s_n)_{n \geq 1}$ are also called ``univoque''
(see \cite{DK1}). Nevertheless, for simplicity we keep our notion of ``univoque".

\begin{definition}\label{selfbrack} 
We define two sets of binary sequences as follows.

\begin{itemize}

\item A sequence $s = (s_n)_{n \geq 1}$ in $\{0,1\}^{{\mathbb N}_+}$ 
    is called {\em self-bracketed} if for every $k \ge 1$
$$
\bar{s} \le_{lex} \sigma^k(s) \le_{lex} s
$$

The set of self-bracketed sequences in $\{0,1\}^{{\mathbb N}_+}$ is 
    denoted by $\Gamma$.

\item If all the inequalities above are strict, the sequence $s$ is said
to be {\em strictly self-bracketed}. If one of the inequalities is an equality,
then $s$ is said to be {\em periodic self-bracketed}.

The subset of $\Gamma$ consisting of the strictly self-bracketed sequences is
denoted by $\Gamma_{strict}$.

\end{itemize}
\end{definition}

\begin{remark}\label{periodic}
The reader will have noticed that the expression ``periodic self-bracketed''
comes from the fact that $\sigma^k(s) = s$ or $\sigma^k(s) = \bar{s}$ for
some $k \geq 1$ implies that the sequence $s$ is periodic.
\end{remark}

With this terminology we can rephrase the following result from \cite{EJK}.

\begin{theorem}\label{EJK}
{\rm \cite{EJK}} A sequence in $\{0,1\}^{{\mathbb N}_+}$ is the unique
$\beta$-expansion of $1$ for a univoque number $\beta$ in $(1,2)$ if and only 
if it strictly self-bracketed.
\end{theorem}

\begin{corollary}\label{cor:01}
Let $s = (s_n)_{n\geq 1}$ be a sequence in $\{0,1\}^{\BbN_{+}}$.
Suppose that the largest string of consecutive $1$'s in $s$ has length $k$,
    and the largest string of consecutive $0$'s has length $n$
(here $k$ and $n$ may be $\infty$.)
If $n > k$, then $s$ is not self-bracketed.
\end{corollary}

\medskip

There exists a smallest univoque real number in $(1, 2)$,
\cite{KL}. Recall first that the Thue-Morse sequence is the fixed point 
beginning with 0 of the morphism $0 \to 01$, $1 \to 10$ (see for example 
\cite{AS}), hence the sequence 
$$
0 \ 1 \ 1 \ 0 \ 1 \ 0 \ 0 \ 1 \ 1 \ 0 \ 0 \ 1 \ 0 \ 1 \ 1 \ 0 \ \ldots
$$
\begin{theorem}\label{KL}
{\rm \cite{KL}} There exists a smallest univoque real number $\kappa \in (1,2)$. 
One has $\kappa \approx 1.787231$, and 
$d_\kappa(1)=(t_n)_{n \ge 1}$, where $(t_n)_{n \ge 1}=11010011\ldots$ is 
obtained by shifting the Thue-Morse sequence.
\end{theorem}
The number $\kappa$ is not rational; actually more can be proved.

\begin{theorem}\label{AC}
{\rm \cite{AC2}} The Komornik-Loreti constant $\kappa$ is transcendental.
\end{theorem}

\noindent
\subsection*{Notation}

\smallskip

In the remainder of this paper, we will denote 
    by ${\mathcal U}$ the set of real numbers
in $(1, 2)$ which are univoque. We will denote by $\widetilde{\mathcal U}$
the set of real numbers $\beta \in (1, 2)$ such that $\ddb$ is finite and $\ddbs$
is a periodic self-bracketed sequence.

Formally, we have 
\[ {\mathcal U} = \{\beta \in (1,2): d_\beta(1) \in \Gamma_{strict}\} \]
and 
\[ \widetilde{\mathcal U} = 
    \{\beta \in (1,2): d_\beta(1)\ \mathrm{is\ finite\ and}\ d_{\beta}^*(1)\ 
            \mathrm{is\ periodic\ self-bracketed}\} \]

\subsection{Pisot numbers}

A {\em Pisot number} is an algebraic integer $>1$ such that all its algebraic 
conjugates (other than itself) have modulus $<1$. As usual the set of Pisot 
numbers is denoted by $S$ and its derived set (set of limit points) by $S'$.
It is known that $S$ is closed \cite{S}, and has a smallest element,
which is the root $>1$ of the polynomial $x^3-x-1$ (approx.~1.3247).
A {\em Salem number} is an algebraic integer $>1$ such that 
all its algebraic conjugates have modulus $\le 1$, with at least one conjugate
on the unit circle. 

We recall some results on Pisot and Salem numbers (the 
reader is referred to~\cite{Ber} for more on these topics).
One important result is that if $\beta$ is a Pisot number then $\ddb$ is 
eventually periodic (finite or infinite) \cite{Be1}. 
Note that $\ddb$ is never periodic, but that when $\ddb$ is finite,
$\ddbs$ is periodic. A number $\beta$ such that $\ddb$ is 
eventually periodic is called a {\em Parry number} (they are called
{\em beta-numbers} by Parry~\cite{Pa}).
When $\ddb$ is finite, $\beta$ is called a {\em simple} Parry number. 

One deeper result is the following one. 

\begin{theorem}\label{AB}
{\rm \cite{Be1, schmidt}} Let $\beta$ be a Pisot number. A number $x$ of $[0,1]$
has a (finite or infinite) eventually periodic greedy $\beta$-expansion
if and only if it belongs to $\Q(\beta)$.
\end{theorem}

For lazy expansions we have a similar result.

\begin{corollary}\label{lep}
Let $\beta$ be a Pisot number. A number $x$ of $[0,1]$
has an eventually periodic lazy $\beta$-expansion
if and only if it belongs to $\Q(\beta)$.
\end{corollary}
\begin{proof}
Let $\ell_\beta(x)=(x_n)_{n\ge 1}$. By Theorem~\ref{EJK1} the sequence
$(\overline{x_n})_{n\ge 1}$ is the greedy $\beta$-expansion of the number
$\frac{\lfloor \beta \rfloor}{\beta -1}-x$, and the result follows from 
Theorem~\ref{AB}.
\end{proof}

\bigskip

Amara has determined all the limit points of $S$ smaller than 2 in~\cite{Am}.
\begin{theorem}\label{Am} {\rm \cite{Am}}
The limit points of $S$ in $(1,2)$ are the following:
$$\varphi_1=\psi_1<\varphi_2<\psi_2<\varphi_3<\chi
<\psi_3<\varphi_4< \cdots <\psi_r<\varphi_{r+1}< \cdots <2$$
where 
$$
\begin{cases}
\text{the minimal polynomial of\ } \varphi_r \text{\ is\ } x^{r+1}-2x^r+x-1, \\
\text{the minimal polynomial of\ } \psi_r \text{\ is\ } x^{r+1}-x^r-\cdots-x-1, \\
\text{the minimal polynomial of\ } \chi \text{\ is\ } x^4-x^3-2x^2+1. \\
\end{cases}
$$
\end{theorem}

The first few limit points are:
\begin{itemize}
\item $\varphi_1 = \psi_1 \approx 1.618033989$, the root in $(1,2)$ of $x^2-x-1$ 
\item $\varphi_2 \approx 1.754877666$, the root in $(1,2)$ of $x^3-2 x^2+x-1$
\item $\psi_2 \approx 1.839286755$, the root in $(1,2)$ of $x^3-x^2-x-1$ 
\item $\varphi_3 \approx 1.866760399$, the root in $(1,2)$ of $x^4-2 x^3+x-1$ 
\item $\chi \approx 1.905166168$, the root in $(1,2)$ of $x^4-x^3-2 x^2+1$ 
\item $\psi_3 \approx 1.927561975$, the root in $(1,2)$ of $x^4-x^3-x^2-x-1$ 
\end{itemize}

The greedy and lazy $\beta$-expansions of these points 
are given in 
    Table \ref{tab:S limit} below.

For any interval $[a,b]$, with $b<2$, an algorithm of Boyd 
\cite{Boyd78, Boyd83, Boyd84} finds all Pisot numbers in the interval.
If $[a,b]$ contains a limit point $\theta$, then there exists an
$\varepsilon >0$ such that all Pisot numbers in
$[\theta-\varepsilon, \theta+\varepsilon]$ are {\em regular\,} Pisot numbers 
of a known form. Boyd's algorithm detects these regular Pisot numbers, 
and truncates the search accordingly. (For a non-effective study of
Pisot numbers in subintervals of $(1,2)$, see also \cite{Ta1, Ta2}.)

   \bigskip
    
Recall that Boyd has shown that for any Salem number of degree 4 the greedy
expansion of 1 is eventually periodic, \cite{Boyd89}, and has given some evidence 
in favor of the conjecture that it is still the case for degree 6, \cite{Boyd96a}.

\section{Preliminary combinatorial results}

We start by defining a function $\Phi$ on the infinite words of the form 
$b=(z0)^\infty$.

\begin{definition}\label{Phi}
Let $b=(z0)^\infty$ be a periodic binary word whose period pattern ends in 
a $0$.  Suppose furthermore that the minimal period of $b$ is equal to $1+|z|$.
Then we define $\Phi(b)$ by
$$
\Phi(b):=(z1\overline{z}0)^\infty.
$$
\end{definition}

We now recall a result from \cite{All}.

\begin{lemma}\label{all}

\mbox{ }

\begin{itemize}

\item If a sequence $b$ belonging to $\Gamma$ begins with 
    $u \bar{u}$ where $u$ is
a finite nonempty word, then $b = (u \bar{u})^\infty$.

\item If $b=(z0)^\infty$, where the minimal period of $b$ is equal to $1+|z|$,
is an element of $\Gamma$, then $\Phi(b)$ belongs to $\Gamma$, and there is
no element of\, $\Gamma$ lexicographically between $b$ and $\Phi(b)$.

\end{itemize}
\end{lemma}

\begin{remark}\label{closed}
The inequalities defining the set $\Gamma$ show that $\Gamma$ is a 
(topologically) closed set.
\end{remark}

\begin{corollary}\label{all2}
Let $b=(z0)^\infty$ (where the minimal period of $b$ is equal to $1+|z|$).
The sequence $(\Phi^{(n)}(b))_{n \geq 0}$ is a sequence of elements of
$\Gamma$ that converges to a limit $\Phi^{(\infty)}(b)$ in $\Gamma$.
The only elements of $\Gamma$ lexicographically between $b$ and 
$\Phi^{(\infty)}(b)$ are the $\Phi^{(k)}(b)$, $k \geq 0$.
\end{corollary}
By abuse of notation, if $\theta$ is the number such that $d^*_\theta(1)=b$,
we denote by $\Phi(\theta)$ the real number $> 1$ such that
$d^*_{\Phi(\theta)}(1)=\Phi(b)$.

\bigskip

Take $b=d^{*}_{\psi_r}(1)=(1^{r}0)^\infty$.
Then $\Phi(b)=(1^{r}10^r0)^\infty = d^{*}_{\varphi_{r+1}}(1)$, thus
$\varphi_{r+1}=\Phi(\psi_r)$. Let $\pi_{r}$ be the real number defined by
$d^{*}_{\pi_{r}}(1)=\Phi^{(\infty)}((1^{r}0)^\infty)$, that is, $\pi_{r} =
    \Phi^{\infty}(\psi_r)$. Then $d^{*}_{\pi_{r}}(1)$
is strictly self-bracketed (see \cite{All}), hence the following result holds true.
\begin{proposition}\label{pi}
The number $\pi_{r}$ is univoque.
Furthermore 
    between $\psi_r$ and $\pi_{r}=\Phi^{(\infty)}(\psi_r)$ the 
    only real numbers
belonging to ${\mathcal U}$ or $\widetilde{\mathcal U}$ are the 
    numbers $\varphi_{r+1}$,
$\Phi(\varphi_{r+1})$, $\Phi^{(2)}(\varphi_{r+1})$, etc. They all belong to
$\widetilde{\mathcal U}$.
\end{proposition}

\bigskip

We will now prove a combinatorial property of the sequences $d_{\beta}(1)$.
Before stating and proving this property we first make a straightforward
remark.

\begin{remark}\label{real}
Let $u$ and $v$ be two binary words having the same length. Let $a$ and $b$
be either two binary words having the same length or two infinite binary 
    sequences.
Suppose that $a$ begins with $u$ and $b$ begins with $v$.
Then
$$
\begin{array}{lll}
a \leq_{lex} b & \Longrightarrow &  u \leq_{lex} v \\
u <_{lex} v    & \Longrightarrow &  a <_{lex} b.
\end{array}
$$
\end{remark}

\begin{proposition}\label{pral}
Let $a=(w0)^{\infty}$ be an infinite periodic binary sequence with
minimal period $1+|w|$, such that $w$ (and hence $a$) 
    begins in $1$. 
Let $b=w10^{\infty}$.
Then the following two properties are equivalent:

(i) $\forall k \geq 1$,\ \ $\sigma^k(a)\leq_{lex} a$,

(ii) $\forall k \geq 1$,\ \ $\sigma^k(b)<_{lex} b$.
\end{proposition}

\begin{proof}

We first prove (i) $\Longrightarrow$ (ii).
Since we clearly have $\sigma^k(b)<_{lex} b$ for each $k \geq |w|$, we can
suppose that $k < |w|$. 
Write $w = u v$ where $|u| = k$, hence $u$ and $v$ are both nonempty. 
This gives $a = (uv0)^{\infty}$ and $b = uv10^{\infty}$,
and we want to prove that $v10^{\infty} <_{lex} uv10^{\infty}$.

Let us write $|v| = d|u| + e$, where $d \geq 0$ and $e \in [0, |u|)$.
We can write
$v = v_1 v_2 \ldots v_d z$, with $|v_1| = |v_2|= \ldots = |v_d| = |u|$,
and $|z| = e < |u|$.
Note that, if $d=0$, then $v=z$.

Let us also write $u=st$ and, for each $j \in [1,d]$, $v_j = s_jt_j$,
where $|s| = |s_1| = |s_2|= \ldots = |s_d| = |z|$ and
$|t| = |t_1| = |t_2| = \ldots = |t_d|$.
We thus have
$$
a = (s t s_1 t_1 s_2 t_2 \ldots s_d t_d z 0)^{\infty}
$$
and we want to prove that
$$
s_1 t_1 s_2 t_2 \ldots s_d t_d z 1 0^{\infty} <_{lex}
s t s_1 t_1 s_2 t_2 \ldots s_d t_d z 1 0^{\infty}.
$$
Applying, for each $j \in [1,d]$, the hypothesis $\sigma^k(a)\leq_{lex} a$
with $k = |s t s_1 t_1 s_2 t_2 \ldots s_{j-1} t_{j-1}|$ (in particular
if $j=1$ then $k = |st|$), we see that $s_j t_j \leq_{lex} s t$. Define
$$
{\mathcal E} := \{j, \ s_j t_j <_{lex} s t\}.
$$

\begin{itemize}

\item If ${\mathcal E} \neq \emptyset$, let $j_0 = \min {\mathcal E}$. Then
$$
s t = s_1 t_1 = s_ 2 t_2 = \ldots = s_{j_0-1} t_{j_0-1}
$$
i.e.,
$$
s = s_1 = s_2 = \ldots = s_{j_0-1} \ \mbox{\rm and } \
t = t_1 = t_2 = \ldots = t_{j_0-1}
$$
(this condition is empty if $j_0=1$) and
$$
s_{j_0} t_{j_0} <_{lex} s t.
$$
In this case we have $b = (st)^{j_0}s_{j_0}t_{j_0} \ldots s_d t_d z10^{\infty}$
and we want to prove that
$$
(st)^{j_0-1}s_{j_0}t_{j_0} \ldots s_d t_d z10^{\infty}
<_{lex} (st)^{j_0}s_{j_0}t_{j_0} \ldots s_d t_d z10^{\infty}
$$
which is an immediate consequence of the inequality $s_{j_0}t_{j_0} <_{lex} st$.

\item If ${\mathcal E} = \emptyset$, then either $d=0$, or
$s_1 t_1 = s_2 t_2 = \ldots = s_d t_d = s t$.
Either way, we get 
$$
s_1 = s_2 = \ldots = s_d = s \ \mbox{\rm and }\
t_1 = t_2 = \ldots = t_d = t.
$$
In this case we have $a = ((st)^{d+1}z0)^{\infty}$ and we want to prove
that $(st)^d z10^{\infty} <_{lex} (st)^{d+1} z10^{\infty}$, i.e., that
$z10^{\infty} <_{lex} st z10^{\infty}$. Applying the hypothesis
$\sigma^k(a)\leq_{lex} a$ with $k = |(st)^{d+1}|$, we see that $z \leq_{lex} s$.

\begin{itemize}

\item If $z <_{lex} s$, the inequality $z10^{\infty} <_{lex} st z10^{\infty}$ is 
clear.

\item If $z = s$, we want to prove that $10^{\infty} <_{lex} t z10^{\infty}$,
i.e., that $t$ begins in $1$ (note that, if $t$ is empty, then the
inequality is clear since $z=s$ begins in $1$ as does $a$). If we had $t=0r$,
with $r$ possibly empty, we would have $a = ((z0r)^{d+1}z0)^{\infty}$.
Applying the hypothesis $\sigma^k(a)\leq_{lex} a$ with $k = |(z0r)^{d+1}|$
and $k = |(z0r)^dz0|$ we get respectively $z0z0r \leq_{lex} z0rz0$
(i.e., $z0r \leq_{lex} rz0$) and $rz0 \leq_{lex} z0r$. Hence we have $rz0 = z0r$.

Writing this last equality as $r (z0) = (z0) r$, the Lyndon-Sch\"utzenberger
theorem (see~\cite{LS62}) implies that $r = \emptyset$ or there exist a nonempty
word $x$ and two integers $p, q \geq 1$, such that $z0 = x^p$ and $r = x^q$.
This gives $a = (x^{a(d+2)})^{\infty}$ or $a = (x^{(p+q)(d+1)+p})^{\infty}$.
In both cases $a = x^{\infty}$ and $|x| < |((z0r)^{d+1}z0)|$ which
contradicts the minimality of the period of $a$.

\end{itemize}

\end{itemize}

We now prove (ii) $\Longrightarrow$ (i).
Because of the periodicity of the sequence $a$ and the fact that it begins in $1$,
we can suppose that $k \leq |w|$. Hence we write $w = uv$ with $u$ and  $v$
nonempty and $|u| = k$, and we want to prove that 
$v0(uv0)^{\infty} \leq_{lex} (uv0)^{\infty}$.
Since $u$ begins in $1$ as $a$ does, it suffices to prove that
$v01^{\infty} \leq_{lex} (uv0)^{\infty}$. Applying the hypothesis
$\sigma^k(b)<_{lex} b$ with $k = |u|$, we have
$v10^{\infty} <_{lex} uv10^{\infty}$.

Hence $v10^{|u|} \leq_{lex} uv1$. This inequality must be strict since its
left-hand side ends in a $0$ and its right-hand side ends with a $1$: thus
$v10^{|u|} <_{lex} uv1$. Hence $v10^{|u|} \leq_{lex} uv0$.

We then can write $v0u <_{lex} v10^{|u|} \leq_{lex} uv0$, hence
$v0u <_{lex} uv0$.  This implies in turn
$v0(uv0)^{\infty} \leq_{lex} (uv0)^{\infty}$.
\end{proof}
\begin{corollary}\label{dstar}
The sequence $a=(w0)^\infty$ is equal to $d_\theta^*(1)$
for some $\theta>1$ if and only if $b=w10^\infty$ is equal to
$d_\theta(1)$.
\end{corollary}

We end this section with a result on limits of sequences of elements in
$\Gamma$.

\begin{lemma}\label{lsup}
A sequence of $\Gamma$ of the form $(w0)^\infty$ cannot be a limit from 
above of a non-eventually constant sequence of elements of $\Gamma$.
\end{lemma}

\begin{proof}
Suppose we have a sequence $(z^{(m)})_{m \ge 0}$ with 
    $z^{(m)} = (z^{(m)}_n)_{n \geq 1}$ belonging to $\Gamma$, 
    and converging towards $(w0)^\infty$,
    with $z^{(m)} \ge (w0)^\infty$.
 From Lemma~\ref{all} there is no element of $\Gamma$ lexicographically
between $(w0)^\infty$ and $(w1\bar{w}0)^\infty$, hence $(z^{(m)})_{m \ge 0}$ 
is ultimately equal to $(w0)^\infty$.
\end{proof}

\section{First results}

In this section we consider only numbers $\beta$ belonging to $(1, 2)$.

\subsection{Preliminary Results}

Our goal here is to present some simple preliminary data.
In particular, in Table \ref{tab:S limit}, we give the expansions 
    for Pisot numbers in $S'\cap(1,2)$, 
    in Table \ref{tab:PisotUnivoque} we give Pisot numbers of small degree 
    in the interval (1,2), 
    and in Table \ref{tab:SalemUnivoque} we examine Salem numbers of small 
    degree in the interval $(1,2)$.
Some observations that are worth making, based on these tables, include:
\begin{remark} {\ }
\begin{itemize}
\item The golden ratio $\varphi_1=\psi_1$ is the smallest element of\, 
    $\widetilde{\mathcal U}$.
    (This comes straight from Definition \ref{selfbrack}.)
\item There is no univoque Pisot number of degree $2$ or $3$.
\item The number $\chi$ is the unique Pisot number of degree 
$4$ which is univoque.
\item For Pisot numbers $\psi_r$, the lazy expansion coincides with $d_{\psi_r}^*(1)$.
\item There exists a unique Salem number of degree $4$ which
is univoque.
\item Salem numbers greater than the Komornik-Loreti constant
$\kappa$ appear to be univoque (for degrees 4 and 6).
\end{itemize}
\end{remark}

\begin{table}[h]
\begin{tabular}{|lllll|}
\hline
Minimal    & Pisot  & Greedy     & Lazy       & Comment  \\
Polynomial & Number &  expansion &  expansion &   \\
\hline
$x^{r+1} - 2 x^{r} + x - 1$ & $\varphi_r$  & $1^r0^{r-1}1$ & 
                            $1^{r-1}01^\infty$ & 
periodic self-bracketed \\
$x^{r+1} - x^{r} - \cdots - 1$ & $\psi_r$  & $1^{r+1}$ & $(1^r0)^\infty$& 
                            periodic self-bracketed  \\
$x^4-x^3-2 x^2 + 1$ & $\chi$  & $11(10)^\infty$ & $11(10)^\infty$ & univoque \\
\hline
\end{tabular}
\caption{Greedy and lazy $\beta$-expansions of real numbers in $S'\cap(1,2)$.}
\label{tab:S limit}
\end{table}

We also observe the following lemma which is straightforward.
\begin{lemma}\label{unit}
A Parry number which is univoque must be a unit (i.e., an algebraic
integer whose minimal polynomial has its constant term equal to $\pm 1$).
\end{lemma}

For each Pisot or Salem number of degree less than 4 or 6 respectively, we simply compute the 
greedy and lazy expansion, and then compare them to see when they are equal.
To find the list of Pisot numbers, we use the algorithm of Boyd \cite{Boyd78}.
 Although there is no nice algorithm to find Salem numbers in $(1,2)$ of 
fixed degree, for low degree we can use brute force.
Namely, if $P(x) = x^n + a_1 x^{n-1} + \cdots + a_1 x + 1$ is a 
    Salem polynomial with root in (1,2) and
    $Q(x) = x^n + b_1 x^{n-1} + \cdots + b_1 x + 1 = (x+2)(x+1/2)(x+1)^{n-2}$,
    then we have $|a_i| \leq b_i$.
See \cite{Borwein02} for more on bounds of coefficients.

\begin{table}
\begin{center}
\begin{tabular}{|lllll|}
\hline
Minimal polynomial & Pisot number & Greedy & Lazy  & Comment\\
                    &              & expansion & expansion & \\
\hline
$x^2-x-1        $ & 1.618033989 & $11           $ & $01^\infty$ & periodic self-bracketed\\
& & & & \\
$x^3-x-1        $ & 1.324717957 & $10001        $ & $00001^\infty$&       \\
$x^3-x^2-1      $ & 1.465571232 & $101          $ & $001^\infty$ &\\
$x^3-2 x^2+x-1  $ & 1.754877666 & $1101         $ & $101^\infty$ & periodic self-bracketed\\
$x^3-x^2-x-1    $ & 1.839286755 & $111          $ & $(110)^\infty$& periodic self-bracketed \\
& & & & \\
$x^4-x^3-1      $ & 1.380277569 & $1001         $ & $0001^\infty$ &\\
$x^4-2 x^3+x-1  $ & 1.866760399 & $111001       $ & $1101^\infty$ & periodic self-bracketed\\
$x^4-x^3-2 x^2+1$ & 1.905166168 & $11(10)^\infty$ & $11(10)^\infty$ & univoque\\
$x^4-x^3-x^2-x-1$ & 1.927561975 & $1111         $ & $(1110)^\infty$ &  periodic self-bracketed\\
\hline
\end{tabular}
\caption{Greedy and lazy expansions of degree 2, 3 and 4 Pisot numbers.}
\label{tab:PisotUnivoque}
\end{center}
\end{table}

\begin{table}
\begin{center}
\begin{tabular}{|lll{p}{1 in}l|}
\hline
Minimal polynomial & Salem        & Greedy & Lazy  & Comment\\
                    & number       & expansion & expansion & \\
\hline
$x^4-x^3-x^2-x+1$ & 1.722083806 & $1(100)^\infty$ & $101(110)^\infty$ & \\
$x^4-2 x^3+x^2-2 x+1$ & 1.883203506 & $1(1100)^\infty$ & $1(1100)^\infty$ & 
                       univoque  \\ 
& & & & \\

$x^6-x^4-x^3-x^2+1$ & 1.401268368 &
 $1(0010000)^\infty$ & $0010111(1111110)^\infty$ & \\

$x^6-x^5-x^3-x+1$ &  1.506135680 &
    $1(01000)^\infty$ & $01011(11110)^\infty$ & \\

$x^6-x^5-x^4+x^3-x^2-x+1$ & 1.556030191 &
    $1(01001001000)^\infty$ &  
$    
0
1^3 
(01)^2 
(
1^7 
0
1^6 
w 
1^3 
w 
1^6 
0
)^\infty$ & \\

$x^6-x^4-2 x^3-x^2+1$ & 1.582347184 &
     $1(0101000)^\infty$ & $011(110)^\infty$ & \\

$x^6-2 x^5+2 x^4-3 x^3$ 
  & 1.635573130 & $1(1000000100)^\infty$ & $1010101(1101111110)^\infty$ & \\
$\ \ \ \ +2 x^2-2 x+1$  & & & & \\

$x^6-x^5-x^4-x^2-x+1$ & 1.781643599 &
  $1(10100)^\infty$ & $11001(11110)^\infty$ &  \\

$x^6-2 x^5+x^3-2 x+1$ & 1.831075825 &
    $1(10110100)^\infty$ & $1(10110100)^\infty$ & univoque \\

$x^6-x^5-x^4-x^3-x^2-x+1$& 1.946856268 &
 $1(11100)^\infty$ & $1(11100)^\infty$ & univoque \\

$x^6-2 x^5-x^4+3 x^3-x^2$ & 1.963553039 &
 $1(111011100)^\infty$ & $1(111011100)^\infty$ & univoque \\
$\ \ \ \ -2 x+1$ & & & & \\

$x^6-2 x^5+x^4-2 x^3+x^2$ & 1.974818708 &
 $1(111100)^\infty$ & $1(111100)^\infty$ & univoque \\
$\ \ \ \ -2 x+1$ & & & & \\

$x^6-2 x^4-3 x^3-2 x^2+1$ & 1.987793167 &
    $1(1111100)^{\infty}$ & $1(1111100)^{\infty}$ & univoque \\
\hline
\end{tabular}
\caption{Greedy and lazy expansions of degree $4$ and $6$ Salem numbers. \newline
         Here $w = 0 11 0 1 0 11 0$.}
\label{tab:SalemUnivoque}
\end{center}
\end{table}

\subsection{Limit points of univoque numbers}

In this section we concern ourselves with the structure of ${\mathcal U} \cap S$
    and $\widetilde{\mathcal U} \cap S$, as well as intersections with 
    the derived set $S'$.
We begin with the following result.

\begin{proposition}\label{lim}
The limit of a sequence of real numbers belonging to ${\mathcal U}$ belongs 
to ${\mathcal U}$ or $\widetilde{\mathcal U}$.
\end{proposition}

\begin{proof}
Let $(\theta_j)_{j \ge 1}$ be a sequence of numbers belonging to ${\mathcal U}$
such that $\lim_{j \to \infty} \theta_j = \theta$. 
Let $a^{(j)}=(a^{(j)}_n)_{n \ge 1} := d_{\theta_j}(1)$. Up to replacing the
sequence $(\theta_j)_{j \ge 1}$ by a subsequence, we may assume that the sequence
of sequences $(a^{(j)}_n)_{n \ge 1}$ converges to a limit $a = (a_n)_{n \ge 1}$ 
when $j$ goes to infinity. Then (dominated convergence):
$$
1=\sum_{n \ge 1} \frac{a_n}{\theta^n}.
$$
For every $j \ge 1$ the number $\theta_j$ belongs to ${\mathcal U}$. Hence 
the sequence $a^{(j)}$ belongs to $\Gamma_{strict}$ hence to $\Gamma$. 
Thus the limit $a=\lim_{j \to\infty} a^{(j)}$ belongs to $\Gamma$ 
(see Remark~\ref{closed}).

If $a$ belongs to $\Gamma_{strict}$, then it is the $\theta$-expansion 
of $1$, and $\theta$ belongs to ${\mathcal U}$.

If $a$ is periodic self-bracketed, it is of the form $a=(w0)^\infty$, where we 
may assume that the minimal period of $a$ is $1+|w|$. From Corollary~\ref{dstar},
$a=d_\theta^*(1)$, $b:=w10^\infty=d_\theta(1)$, and $\theta$ belongs to the set
$\widetilde{\mathcal U}$.
\end{proof}

\begin{corollary}\label{co}
The numbers $\varphi_{r}$ cannot be limit points of numbers in ${\mathcal U}$.
\end{corollary}

\begin{proof}
This is a consequence of the first part of Lemma~\ref{all}. \end{proof}

\bigskip

We now give two remarkable sequences of real numbers that converge 
to the Komornik-Loreti constant $\kappa$. Part (ii) of 
Proposition~\ref{tr} below was obtained independently by the 
second author and in~\cite{KLP}.

\begin{proposition}\label{tr}

\mbox{ }

\begin{itemize}

\item[(i)] Let $t=(t_n)_{n \ge 1}=11010011\ldots$ be the shifted 
Thue-Morse sequence, and let $\tau_{2^k}$ be the real number $>1$ 
such that $d_{\tau_{2^k}}(1)=t_1 \cdots t_{2^k}$. Then, the sequence 
of real numbers $(\tau_{2^k})_{k \ge 1}$ converges from below to the 
Komornik-Loreti constant $\kappa$.
These numbers belong to $\widetilde{\mathcal U}$. The first three are Pisot 
numbers.

\item[(ii)] There exists a sequence of univoque Parry numbers
that converges to $\kappa$ from above.

\end{itemize}

\end{proposition}

\begin{proof}

To prove (i), note that $\tau_2$ is the golden ratio, 
$\tau_4=\Phi(\tau_2)=\varphi_2$, $\tau_8=\Phi^2(\tau_2)$, etc., 
and $\kappa=\Phi^{(\infty)}(\tau_2)$. 

\bigskip

In order to prove (ii) we define $\delta_{2^k}$ as the number such that
$$
d_{\delta_{2^k}}(1)=
t_1 \cdots t_{2^k-1}(1\overline{t_1} \cdots \overline{t_{2^k-1}})^\infty.
$$
Clearly the sequence $d_{\delta_{2^k}}(1)$ converges to $t$ when $k$ goes
to infinity and thus the sequence $(\delta_{2^k})_{k \ge 1}$ converges to
$\kappa$. \end{proof}

\begin{remark}\label{ass1-ass2}

\mbox{ }

\begin{itemize}
\item Let $Q_{2^k}$ be the polynomial ``associated'' with $\tau_{2^k}$: 
writing $1=\sum_{1 \leq j \leq 2^k} \dfrac{t_j}{\tau_{2^k}^j}$ immediately gives a polynomial
$Q_{2^k}(x)$ of degree $2^k$ such that $Q_{2^k}(\tau_{2^k}) = 0$. Then,
for $k \ge 2$, the polynomial $Q_{2^k}(x)$ is divisible by the product
$(x+1)(x^2+1) \cdots (x^{2^{k-2}}+1)$.
\item Let $R_{2^k}(x)$ be the polynomial of degree $2^{k+1}-1$ 
associated (as above) with $\delta_{2^k}$. Then it can be shown that, 
for $k \ge 2$, the polynomial $R_{2^k}(x)$ is divisible by the same product $
(x+1)(x^2+1) \cdots (x^{2^{k-2}}+1)$.
\end{itemize}

\end{remark}
 
\section{Main results}

Recall that Amara gave in \cite{Am} a complete description of the limit 
points of the Pisot numbers in the interval $(1,2)$ (see Theorem~\ref{Am}).
Talmoudi \cite{Ta2} gave a description for sequences of Pisot numbers approaching
each of the values $\varphi_r, \psi_r$ or $\chi$. The Pisot numbers in these 
sequences are called {\em regular Pisot numbers}.
Further, Talmoudi showed that, for all $\varepsilon > 0$, there are only a 
finite number of Pisot numbers in $(1, 2-\varepsilon)$, that are not in one 
of these sequences. These are called the {\em irregular Pisot numbers}, and 
they will be examined later in Section~\ref{sec:finite}.

Since $\chi$ is a univoque Pisot number (Tables \ref{tab:S limit} and 
    \ref{tab:PisotUnivoque}),
it is natural to ask if there are any other univoque Pisot numbers
smaller than $\chi$. As well, it is natural to ask if there is a smallest 
univoque Pisot number. This leads us to our first result:

\begin{theorem}\label{uni}
There exists a smallest Pisot number in the set\, ${\mathcal U}$.
\end{theorem}

\begin{proof}
Define $\theta$ by $\theta := \inf(S \cap {\mathcal U})$. We already know  
that $\theta$ belongs to $S$, since $S$ is closed. On the other hand, from
Proposition~\ref{lim}, either $\theta$ belongs to ${\mathcal U}$ or to 
$\widetilde{\mathcal U}$. It suffices to show that $\theta$ cannot belong to
$\widetilde{\mathcal U}$. If it were the case, first $\theta$ would be a limit 
point of elements of $(S \cap {\mathcal U})$. On the other hand we could write 
$d^*_{\theta}(1) = (w0)^\infty$, with the minimal period of the sequence 
$d^*_{\theta}(1)$ being $1+|w|$ (note that $\theta < \chi$ since $\chi$ 
belongs to $(S \cap {\mathcal U})$ and $\theta \neq \chi$). 
But from Lemma~\ref{lsup} there is a contradiction.
\end{proof}

\bigskip

Now, to find the univoque Pisot numbers less than $\chi$, we need
    to examine the irregular Pisot numbers less than $\chi$
    (Section~\ref{sec:finite}).
We need also to examine the infinite sequences of Pisot numbers tending 
to those $\varphi_r$ and $\psi_r$ less than $\chi$.
Lastly, we need to examine the sequences of Pisot numbers tending to $\chi$
 from below.

By noticing that $\varphi_1 = \psi_1$ and $\varphi_2$ are all 
    strictly less than $\kappa$, 
    the Komornik-Loreti constant, we can disregard these limit points.
Further, we may disregard $\varphi_3$ as a limit point by Corollary
    \ref{co}.
In particular:
\begin{proposition}\label{thm:gamma}
There are no univoque numbers between $\psi_2$ and $1.8705$. (Note that 
$1.8705 > \varphi_3$.)
\end{proposition}

\begin{proof}
We easily see from Proposition~\ref{pi} 
    that \[ \Phi^{2}(\psi_2) = \Phi(\varphi_3) \approx 1.870556617 \] 
which gives the result.
\end{proof}

\bigskip

So we see that it suffices to examine the sequence of Pisot numbers tending 
towards $\psi_2$ from below, and those tending to $\chi$ from below.

\subsection{Approaching $\psi_2$ from below}
\label{sec:psi}

We know that the $\psi_r$ are limit points of the set of Pisot numbers.
Moreover, we know exactly what the sequences tending to $\psi_r$ look 
like. Let $\Pbr(x) = x^{r+1} - \cdots - 1$ be the Pisot polynomial associated 
with $\psi_r$.
Let $\Abr(x) = x^{r+1} -1$ and $\Bbr(x) = \frac{x^r-1}{x-1}$ be two polynomials 
associated with $\Pbr(x)$.\footnote{Note that the definition of $\Bbr(x)$ is 
different from the definition in \cite{Boyd96}, and corrects a misprint in that paper.}
Then for sufficiently large $n$, 
    the polynomials $\Pbr(x) x^n \pm \Abr(x)$ and 
    $\Pbr(x) x^n \pm \Bbr(x)$
    admit a unique root between $1$ and $2$, which is a Pisot number.
These sequences of Pisot numbers are the regular Pisot numbers 
    associated with $\psi_r$.
See for example \cite{Am, Boyd96}.

Moreover, we have that the roots of
    $\Pbr(x) x^n - \Abr(x)$ and $\Pbr(x) x^n - \Bbr(x)$ approach $\psi_r$ 
    from above, and those of
    $\Pbr(x) x^n + \Abr(x)$ and $\Pbr(x) x^n + \Bbr(x)$ approach $\psi_r$ 
    from below.
This follows as $\Pbr(1) = -1$ and $\Pbr(2) = 1$, with $\Pbr(x)$ strictly
    increasing on $[1,2]$, along with the fact that 
    on $(1, 2]$ we have $\Abr(x), \Bbr(x) > 0$.
Although we need only examine the sequences of Pisot numbers approaching 
    $\psi_2$ from below, we give the results for all sequences 
    approaching $\psi_2$ for completeness.

\begin{lemma}\label{psi_2}
The greedy and lazy expansions of Pisot numbers approaching $\psi_2$ are
    summarized in Table \ref{tab:GL psi_2}.
\end{lemma}

\begin{remark}
It is interesting to observe that, in the case 
    $P_{\psi_2}(x) x^n - B_{\psi_2(x)}(x)$, 
   (last section of Table \ref{tab:GL psi_2}), 
    for $n = 2, 3$ and $4$, 
    the lazy expansion $\ell_\beta(1)$ is equal to $d_\beta^*(1)$.
\end{remark}

\begin{sidewaystable}
\centering
\begin{tabular}{|lllr|}
\hline
Case & Greedy expansion  & Lazy expansion & Comment \\ 
\hline
\multicolumn{4}{|c|}{$P_{\psi_2}(x) x^n + A_{\psi_2}(x)$} \\ &&&\\
$n = 1 $& $101 $ & $00(1)^{\infty} $ & \\
$n = 2 $& $10101 $ & $0(11101)^{\infty} $ & \\
$n = 3 $& $110001 $ & $1010(1)^{\infty} $ & \\
$n = 4 $& $1100110001 $ & $10(1111110011)^{\infty} $ & \\
$n = 3 k + 1 $& $(110)^{k} 011(000)^{k} 1 $ & $(110)^{k} 0(101)^{\infty} $ & \\
$n = 3 k + 2 $& $1(101)^{k} 010(000)^{k} 1 $ & $1(101)^{k} 0(011)^{\infty} $ & \\
$n = 3 k + 3 $& $(110)^{k+1} 00(000)^{k} 1 $ & $1(101)^{k} 0((110)^{k+1} 01(101)^{k} 1)^{\infty} $ & \\
\hline
\multicolumn{4}{|c|}{$P_{\psi_2}(x) x^n - A_{\psi_2}(x)$} \\ &&&\\
$n = 1 $& \multicolumn{3}{l|}{Root bigger than 2} \\
$n = 2 $& \multicolumn{3}{l|}{Root bigger than 2} \\
$n = 3 $& $111(110)^{\infty} $ & $111(110)^{\infty} $ & univoque \\
$n = 4 $& $111(0110)^{\infty} $ & $111(0110)^{\infty} $ & univoque \\
$n = 3 k + 1 $& $111(0(000)^{k-1} 110)^{\infty} $ & $11(011)^{k-1} 1((011)^{k} 0)^{\infty} $ & \\
$n = 3 k + 2 $& $111(00(000)^{k-1} 110)^{\infty} $ & $11(011)^{k-1} 1001(101)^{k-1} 0111(11(011)^{k-1} 110)^{\infty} $ & \\
$n = 3 k + 3 $& $111((000)^{k} 110)^{\infty} $ & $11(011)^{k} 1(110)^{\infty} $ & \\
\hline
\multicolumn{4}{|c|}{$P_{\psi_2}(x) x^n + B_{\psi_2}(x)$} \\ &&&\\
$n = 1 $& $10001 $ & $0000(1)^{\infty} $ & \\
$n = 2 $& $11 $ & $0(1)^{\infty} $ & periodic self-bracketed \\
$n = 3 $& $11001010011 $ & $(1011110)^{\infty} $ & \\
$n = 4 $& $11010011001011 $ & $110100(10111111)^{\infty} $ & \\
$n = 3 k + 1 $& $1(101)^{k} 00(110)^{k} 0(101)^{k} 1 $ & 
  $(1(101)^{k} 00(110)^{k} 0(101)^{k} 0)^{\infty} $ & periodic self-bracketed \\
$n = 3 k + 2 $& $1(101)^{k} 1 $ 
  & $(1(101)^{k} 0)^{\infty} $ & periodic self-bracketed \\
$n = 3 k + 3 $& $(110)^{k+1} 0(101)^{k+1} 001(101)^{k} 1 $ & 
  $(110(110)^{k} 0(101)^{k+1} 001(101)^{k} 0)^{\infty} $ 
  & periodic self-bracketed \\
\hline
\multicolumn{4}{|c|}{$P_{\psi_2}(x) x^n - B_{\psi_2}(x)$} \\ &&&\\
$n = 1 $& \multicolumn{3}{l|}{Root bigger than 2} \\
$n = 2 $& $11111 $ & $(11110)^{\infty} $ & periodic self-bracketed \\
$n = 3 $& $111011 $ & $(111010)^{\infty} $ & periodic self-bracketed \\
$n = 4 $& $1110011 $ & $(1110010)^{\infty} $ & periodic self-bracketed \\
$n = 3 k + 1 $& $11100(000)^{k-1} 11 $ & $
(11(011)^{k-1} 1001(101)^{k-1} 0)^{\infty} $ & \\
$n = 3 k + 2 $& $111(000)^{k} 11 $ & 
$(11(011)^{k} 11(101)^{k} 0)^{\infty} $ & \\
$n = 3 k + 3 $& $1110(000)^{k} 11 $ & 
$(11(011)^{k} (101)^{k+1} 0)^{\infty} $ & \\
\hline
\end{tabular}
\caption{Greedy and lazy expansion for regular Pisot numbers approaching 
    $\psi_2$.}
\label{tab:GL psi_2}
\end{sidewaystable}

\begin{proof}
Table \ref{tab:GL psi_2}, as well as Table \ref{tab:GL chi} later on, 
    are the results of a computation.  
The results themselves are easy to verify, so the main interest is the 
    process that the computer went through, to discover these results.
This is the subject of Section \ref{sec:computer}.
We also list which of these numbers correspond to periodic self-bracketed 
    sequences for completeness.
\end{proof}

\bigskip

This Lemma gives an easy corollary

\begin{corollary}\label{cor:psi}
There exists a neighborhood $[\psi_2-\varepsilon, \psi_2+\varepsilon]$ 
    that contains no univoque numbers.
\end{corollary}

In fact we will see in Section \ref{sec:finite} that this is actually quite a
    large neighborhood.  
This is probably also true for other $\psi_r$, where the neighborhood
    would not be nearly as large.

\subsection{The limit point $\chi$}

We know that $\chi$ is a limit point of the set of Pisot numbers.
Moreover, we know exactly what the sequences tending to $\chi$ 
    look like.
Let $\Px(x) = x^4-x^3-2 x^2+1$ be the Pisot polynomial associated 
    with $\chi$.
Let $\Ax(x) = x^3+x^2-x-1$ and $\Bx(x) = x^4-x^2+1$
    be two polynomials associated with $\Px(x)$.  
Then for sufficiently large $n$, 
    the polynomials $\Px(x) x^n \pm \Ax(x)$ and $\Px(x) x^n \pm \Bx(x)$
    admit a unique root between $1$ and $2$, which is a Pisot number.
See for example \cite{Am, Boyd96}.
    
Moreover, we have that the roots of
    $\Px(x) x^n - \Ax(x)$ and $\Px(x) x^n - \Bx(x)$ approach $\chi$ 
    from above, and those of
    $\Px(x) x^n + \Ax(x)$ and $\Px(x) x^n + \Bx(x)$ approach $\chi$ 
    from below.
This follows as $\Px(1) = -1$ and $\Px(2) = 1$, with $\Px(x)$ strictly
    increasing on $[1,2]$, along with the fact that 
    on $(1, 2]$ we have $\Ax(x), \Bx(x) > 0$.
Although we need only examine the sequences of Pisot numbers approaching 
    $\chi$ from below, we give the results for all sequences 
    approaching $\chi$ for completeness.

\begin{lemma}\label{chi}
The greedy and lazy expansions of Pisot numbers approaching $\chi$ are
    summarized in Table~\ref{tab:GL chi}.
\end{lemma}

\begin{sidewaystable}
\centering
\begin{tabular}{|ll{p}{2.6in}r|}
\hline
Case & Greedy expansion  & Lazy expansion & Comment \\ 
\hline
\multicolumn{4}{|c|}{$\Px(x) x^n + \Ax(x)$} \\ &&&\\
$n = 1$ & $1001001 $ & $00(1111011)^{\infty} $ & \\
$n = 2$ & $11 $ & $0(1)^{\infty} $ & periodic self-bracketed \\
$n = 4$ & $110110101001001011 $ & $110110100(1)^{\infty} $ & periodic self-bracketed \\
$n = 2 k + 1$ & $11(10)^{k-1} 01000(10)^{k-1} 0(00)^{k} 11 $ 
& $11(10)^{k-1} 00(11)^{k} 00(1)^{\infty} $ & \\
$n = 2 k + 2$ & $11(10)^{k-1} 0111000(10)^{k-2} 000010(00)^{k-1} 11 $ 
& $11(10)^{k-1} 01101(11)^{k-1} 00(1)^{\infty} $ & \\
\hline
\multicolumn{4}{|c|}{$\Px(x) x^n - \Ax(x)$} \\ &&&\\
$n = 1$ & \multicolumn{3}{l|}{Root bigger than 2} \\
$n = 2$ & \multicolumn{3}{l|}{Root bigger than 2} \\
$n = 3$ & \multicolumn{3}{l|}{Root bigger than 2} \\
$n = 5$ & $1111(0001100)^{\infty} $ & $1111000101111(0111100)^{\infty} $ & \\
$n = 2 k + 1$ & $111(01)^{k-2} 1(00011(10)^{k-2} 00)^{\infty} $ 
& $111(01)^{k-2} 100011(10)^{k-3} 0111011(1(01)^{k-3} 11111000)^{\infty} $ & \\
$n = 2 k + 2$ & $111(01)^{k-1} 1011((10)^{k-1} 0111(01)^{k-1} 1000)^{\infty} $
& $111(01)^{k-1} 1011((10)^{k-1} 0111(01)^{k-1} 1000)^{\infty} $ & 
univoque\\
\hline
\multicolumn{4}{|c|}{$\Px(x) x^n + \Bx(x)$} \\ &&&\\
$n = 1$ & $10001 $ & $0000(1)^{\infty} $ & \\
$n = 2$ & $101000101 $ & $0(1101)^{\infty} $ & \\
$n = 3$ & $11001 $ & $10(11011)^{\infty} $ & \\
$n = 4$ & $110101(01100110000100)^{\infty} $ & $110(1010110010110111)^{\infty} $ & \\
$n = 5$ & $1110001 $ & $110(1110111)^{\infty} $ & \\
$n = 2 k + 1$ & $11(10)^{k-1} 001 $ & $1110((10)^{k-3} 01111)^{\infty} $ & \\
$n = 2 k + 2$ & $11(10)^{k-1} 0101(1(10)^{k-2} (011)^2(10)^{k-2} 010^4100)^{\infty} $ 
& $1110((10)^{k-2} 01011(10)^{k-2} (011)^2(10)^{k-2} 001^3(10)^{k-2} 01^4)^{\infty} $ & \\
\hline
\multicolumn{4}{|c|}{$\Px(x) x^n - \Bx(x)$} \\ &&&\\
$n = 1$ & \multicolumn{3}{l|}{Root bigger than 2} \\
$n = 2$ & \multicolumn{3}{l|}{Root bigger than 2} \\
$n = 3$ & \multicolumn{3}{l|}{Root bigger than 2} \\
$n = 4$ & $111111000001 $ & $111110(1)^{\infty} $ & periodic self-bracketed \\
$n = 5$ & $1111001111000001 $ & $111100(11101)^{\infty} $ & \\
$n = 2 k + 1$ & $111(01)^{k-2} 101000(10)^{k-3} 011(1(00)^{k-1} 10)^{\infty} $ & $111(01)^{k-2} 100(1(11)^{k-1} 01)^{\infty} $ & \\
$n = 2 k + 2$ & $111(01)^{k-2} 100000(10)^{k-2} 001 $ & $111(01)^{k-1} 110(1)^{\infty} $ & \\
\hline
\end{tabular}
\caption{Greedy and lazy expansion for regular Pisot numbers approaching $\chi$.}
\label{tab:GL chi}
\end{sidewaystable}

Lemma~\ref{chi} above, along with Proposition~\ref{pi} and 
    Corollary \ref{cor:psi} prove the following result:
\begin{theorem}\label{thm:finite}
There are only a finite number of univoque Pisot numbers less than $\chi$.
\end{theorem}

In addition, Lemma~\ref{chi} proves the result

\begin{theorem}\label{thm:chi}
The univoque Pisot number $\chi$ is the smallest limit point of
univoque Pisot numbers. It is a limit point from above
of regular univoque Pisot numbers.
\end{theorem}

\subsection{Univoque Pisot numbers less than $\chi$}
\label{sec:finite}

Our goal in this section is to describe our search for univoque Pisot numbers
    below the first limit point $\chi$.
We know that all univoque Pisot numbers less than $\chi$ are either 
    in the range $[\kappa, \psi_2]$, or in the range $[\pi_2, \chi]$.
Here $\kappa$ is the Komornik-Loreti constant, (approximately $1.787231$),
    $\psi_2$ is approximately $1.839286755$, 
    $\pi_2 > 1.8705$ and
    $\chi$ is approximately $1.905166168$.
We will search for Pisot numbers in the range $[1.78, 1.85]$ and 
    $[1.87, 1.91]$.

To use the algorithm of Boyd \cite{Boyd78}, we need to 
    do an analysis of the limit points in these two ranges.  
In particular, we need to do an analysis of the limit points 
    $\psi_2$ and $\chi$.

We use the notation of \cite{Boyd78}.
Let $P(z)$ be a minimal polynomial of degree $s$ of a Pisot number $\theta$, 
    and $Q(z) = z^s P(1/z)$ be the reciprocal polynomial.
Let $A(z)$ be a second polynomial with integer coefficients, such that
    $|A(z)| \leq |Q(z)|$ for all $|z| = 1$.
Then $f(z) = A(z)/Q(z) = u_0 + u_1 z + u_2 z^2 + \cdots \in \BbZ[[z]]$ 
    is a rational function associated with $\theta$.
The sign of $A(z)$ is chosen in such a way that $u_0 \geq 1$.
Then by Dufresnoy and Pisot \cite{DufresnoyPisot55} we have the following.
\begin{equation}\label{eq:rules}
\begin{array}{rclcl}
1 &\leq &u_0 \\
u_0^2-1 &\leq &u_1 \\
w_n(u_0, \cdots, u_{n-1}) & \leq & u_n &\leq & w^*_n(u_0, \cdots, u_{n-1}) 
\end{array}
\end{equation}
where $w_n$, and $w_n^*$ are defined below.
Let $D_n(z) = -z^n + d_1 z^{n-1} + \cdots + d_n$ and $E_n(z) = -z^n D_n(1/z)$.
Solve for $d_1, \cdots, d_n$ such that 
\[
\frac{D_n(z)}{E_n(z)} = u_0 + u_1 z + \cdots + u_{n-1} z^{n-1} + 
                        w_n(u_0, \cdots, u_{n-1})z^n + \cdots\]
This will completely determine $w_n$.
There are some nice recurrences for $w_n$ and $D_n$, which simplify the 
    computation of $w_n$ \cite{Boyd78}.
We have that $w_n^*$ is computed very similarly, instead considering 
    $D_n^*(z) = z^n + d_1 z^{n-1} + \cdots + d_n$ and $E_n^*(z) = z^n D_n(1/z)$.
Expansions $u_0 + u_1 z + \cdots$ satisfying Equation~(\ref{eq:rules})
    with integer coefficients are in a one-to-one correspondence with 
    Pisot numbers.

Using this notation, Lemma~2 of \cite{Boyd78} becomes:
\begin{lemma}
Let $f = u_0 + u_1 z + u_2 z^2 + \cdots $ be associated with a 
    limit point $\theta$ in $S'$.
Suppose that $w^*_N - w_N \leq 9/4$ for some $N$.
Then for any $n \geq N$, there are exactly two $g$ with expansions 
    beginning with $u_0 + u_1 z + \cdots + u_{n-1} z^{n-1}$.
Moreover, for all $n \geq N$, all $g$ 
    beginning with $u_0 + u_1 z + \cdots + u_{n-1} z^{n-1}$
    are associated with the regular Pisot numbers approaching 
    the limit point $\theta$.
\end{lemma}

So in particular, we need to find the expansion of the limit points 
    around $\psi_2$ and $\chi$.
Consider the following rational functions associated with the limit points 
    $\psi_2$ and $\chi$.
\begin{enumerate}
\item Consider 
    \[-\frac{x+1}{x^3+x^2+x-1} = 1+ 2 x+ 3 x^2+ 6 x^3+ 11 x^4+
     20 x^5+ 37 x^6+ \cdots\]
    the first of the two rational functions associated with the limit 
    point $\psi_2$.
    A quick calculation shows that $w_{24}^* - w_{24} < 9/4$.
\item Consider
    \[ \frac{x^3-1}{x^3+x^2+x-1} = 
       1+x+2 x^2+3 x^3+6 x^4+11 x^5+20 x^6 +\cdots \]
    the second of the two rational functions associated with the limit point
        $\psi_2$.
    A quick calculation shows that $w_{11}^* - w_{11} < 9/4$.
\item Consider
    \[ -\frac{x^3+x^2-x-1}{x^4-2 x^2-x+1} = 
      1+2 x+3 x^2+6 x^3+11 x^4+21 x^5+40 x^6 +\cdots \]
    th first of the two rational functions associated with the limit 
    poit $\chi$.
    A quick calculation shows that $w_{33}^* - w_{33} < 9/4$.
\item onsider
    \[ \frac{x^4-x^2+1}{x^4-2 x^2-x+1} = 
       1+x+2 x^2+4 x^3+8 x^4+15 x^5+29 x^6+\cdots \]
    the second f the two rational functions associated with the limit point
        $\chi$.
    A quick calculation shows that $w_{44}^* - w_{44} < 9/4$.
\end{enumerate}

Using this result, we were able to use Boyd's algorithm for finding 
    Pisot numbers in the two ranges $[1.78, 1.85]$ and $[1.87, 1.91]$,
    (which contain $[\kappa, \psi_2]$ and $[\pi_2, \chi]$), 
    where when we have an expansion that matches one of the four  
    rational functions listed above, we prun that part of the
    search tree, as we would only find regular Pisot numbers of a 
    known form.

There were 227 Pisot numbers in the first range (minus the known regular
    Pisot numbers pruned by the discussion above), and 303 in the 
    second range (similarly pruned).
There were 530 such Pisot numbers in total.

A corollary of this computation worth noting is
\begin{corollary} {\ }
\begin{itemize}
\item The only Pisot numbers in $[\psi_2-10^{-8}, \psi_2+10^{-8}]$ are 
    $\psi_2$ and the regular Pisot numbers associated with $\psi_2$.
\item The only Pisot numbers in $[\chi-10^{-13}, \chi+10^{-3}]$ are 
    $\chi$ and the regular Pisot numbers associated with $\chi$.
\end{itemize}
\end{corollary}

We then checked each of these 530 Pisot numbers to see if they were univoque.
We did this by computing the greedy and lazy $\beta$-expansion of each 
    Pisot number
    and checked if they were equal.
This calculation gave the following theorem:
    
\begin{theorem}\label{thm:finite2}
There are exactly two univoque Pisot numbers less than $\chi$.
They are
\begin{itemize}
\item $1.880000\cdots$ the root in $(1,2)$ of the polynomial
    ${x}^{14}-2{x}^{13}+{x}^{11}-{x}^{10}-{x}^{7}+{x}^{6}
    -{x}^{4}+{x}^{3}-x+1$ with univoque expansion $111001011(1001010)^\infty$.
\item $1.886681\cdots$ the root in $(1,2)$ of the polynomial
    ${x}^{12}-2{x}^{11}+{x}^{10}-2{x}^{9}+{x}^{8}-{x}^{3}+{x}^{2}-x+1$ 
    with univoque expansion $111001101(1100)^\infty$
\end{itemize}
\end{theorem}
 
\section{Regular Pisot numbers associated with $\psi_r$}

The goal of this section is to show that $2$ is the limit point 
    of univoque Pisot numbers.
We will do this by observing that for each $r$, there are regular Pisot 
    numbers between $\psi_r$ and $2$ that are univoque.
We know that the $\psi_r$ are limit points of the set of regular Pisot 
    numbers.
Moreover we know that $\psi_r \to 2$ as $r \to \infty$.
Using the notation of Section~\ref{sec:psi} we define 
    $\Pbr$ and $\Abr$.
We denote the Pisot number associated with the polynomial 
    $\Pbr(x) x^n - \Abr(x)$, as $\psiAm{r}{n}$.

\begin{theorem}
Let $n \geq r + 1$.
Then the greedy expansion of $\psiAm{r}{n}$ is
\begin{equation} 1^{r+1} (0^{n-r-1} 1^{r} 0)^{\infty}. 
\label{eq:general}
\end{equation}
\end{theorem}

\begin{proof}
First we expand this expansion to see that it is equivalent to
\[
\begin{array}{lrcl}
& 1 & = & \frac{1}{x} + \cdots + \frac{1}{x^{r+1}} + 
          (0 + \frac{1}{x^{n+1}} + \cdots + \frac{1}{n+r} + 0) 
          \left(\frac{1}{1-\frac{1}{x^n}}\right) \\
\implies & 1 & = & 
    \frac{\frac{1}{x}-\frac{1}{x^{r+2}}}{1-\frac{1}{x}} + 
    \left(\frac{\frac{1}{x^{n+1}}-\frac{1}{x^{n+r+1}}}{1-\frac{1}{x}}\right)
    \left(\frac{1}{1-\frac{1}{x^n}}\right) \\
\implies & (1-\frac{1}{x})(1-\frac{1}{x^n}) & = & 
    (\frac{1}{x}-\frac{1}{x^{r+2}}) (1-\frac{1}{x^n}) + 
    (\frac{1}{x^{n+1}}-\frac{1}{x^{n+r+1}}) \\
\implies & x^{r+2} (x-1) (x^n-1) & = & 
    x (x^{r+1}-1)(x^n-1) + x^{r+2} - x^2 \\
\implies & 0 & = & 
    x^{n+r+3} - 2 x^{n+r+2} +x^{n+1} - x^{r+3} + x^{r+2} - x^2 - x \\
\implies & 0 & = & 
    x (x-1) (x^n(x^{r+1} - x^r - \cdots - 1) - (x^{r+2}-1)) \\
\implies & 0 & = & 
    x (x-1)(\Pbr(x) x^n - \Abr(x)) 
\end{array}
\]
So we see that this is a valid expansion for this regular Pisot number.

To observe that this is indeed the greedy $\beta$-expansion we observe that the 
    $\beta$-expansion starts with $r+1$ consecutive $1$'s, and all strings 
    of consecutive $1$'s after this are shorter than $r+1$. 
Hence it follows from Theorem \ref{Pa}.
\end{proof}

By Corollary \ref{cor:01} we get the immediate result:
\begin{corollary}
If $n \geq 2 (r+1)$ the regular Pisot number $\psiAm{r}{n}$ is not univoque.
\end{corollary}

So, the main theorem is
\begin{theorem}
Assume $r+1 \leq n < 2 (r+1)$.
Then $\psiAm{r}{n}$ is a univoque Pisot number.
\end{theorem}

\begin{proof}
So it suffices to see that the equation (\ref{eq:general}) is both greedy
   and lazy. 
This follows from Theorems \ref{Pa} and \ref{EJK1}.
\end{proof}

\begin{corollary}
We have that $2$ is a limit point of $S \cap \mathcal U$.
\end{corollary}

\begin{proof}
We see that $\psiAm{r}{n}$ is always greater than $\psi_r$
Further, for $r+1 \leq n \leq 2(r+1)$ we have that $\psiAm{r}{n}$ is 
    less than 2, which follows from
     noticing that 
    $\Pbr(1) 1^n - \Abr(1) =  (1 - 1 - 1 - \cdots - 1) - (1^{r+1} -1) < 0$ and 
    $\Pbr(2) 2^n - \Abr(2) =  
    2^n(2^{r+1} - 2^{r} - \cdots - 1) - (2^{r+1} - 1) = 
    2^n -  2^{r+1} + 1  > 0$.
Further, we see that $\psi_r$ tends to 2.
\end{proof}

\section{Automated conjectures and proofs}
\label{sec:computer}

The results in Tables \ref{tab:GL psi_2} and \ref{tab:GL chi} were
    generated automatically.
This section describes the algorithms that were used
    to do this.

\begin{itemize}
\item {\bf Computing the greedy $\beta$-expansion.}   
    We will explain, given $\beta$ a root of $P_\beta(x)$, how to compute the 
    greedy $\beta$-expansion of $1$, (assuming periodicity).
\item {\bf Computing the lazy $\beta$-expansion.}
    We will explain, given $\beta$ a root of $P_\beta(x)$, how to compute the 
    lazy $\beta$-expansion of $1$, (assuming periodicity).
\item {\bf Creating the conjecture.}
    We will explain how with the greedy or lazy $\beta$-expansion  of $1$ for 
    a sequence of regular Pisot numbers, how to create a conjecture
    of the general pattern of the $\beta$-expansion.
\item {\bf Verifying conjecture.}
    We will explain how a general pattern can be verified to be a valid
    $\beta$-expansion.
\item {\bf Check greedy/lazy/univoque/periodic self-bracketed expansion.}
    We will explain how to check if a general pattern is a valid greedy, lazy,
    univoque or periodic self-bracketed $\beta$-expansion.
\end{itemize}
 
\subsection{Computing the greedy $\beta$-expansion}
\label{sec:Greedy Comp}

The greedy algorithm does the most work 
   possible at any given step (see discussion in Section \ref{sec:greedy}).

The computation is done symbolically modulo the minimal polynomial of $\beta$,
    and floating point numbers are used only when computing $x_n$.
A check is done on $\beta r_{n-1} - x_n$ to ensure that the 
    calculation is being done with sufficient digits to guarantee
    the accuracy of the result.

A list of previously calculated $r_n$'s is kept and checked upon
    each calculation to determine when the $\beta$-expansion becomes
    eventually periodic.

\subsection{Computing the lazy $\beta$-expansion}
\label{sec:Lazy Comp}

Basically, the algorithm tries to do the minimal work at any given time
    (see discussion in Section \ref{sec:lazy}).

As with the greedy expansion,  computations are done as a mixture of floating point and symbolic,
    to allow for recognition of periodicity,
    with the same checks being performed as before 
    to ensure the accuracy of the result.

\subsection{Creating the conjecture}
\label{sec:conj}

In this section we will explain how, given $d_{q_1}(1)$ and $d_{q_2}(1)$, (or 
    the related lazy $\beta$-expansions), for some ``regular sequence'' of Pisot 
    numbers $q_k$, we can conjecture a ``nice'' expression for $d_{q_k}(1)$.
This is probably best done by example.
Assume that two consecutive greedy expansions give the finite expansions:
\begin{eqnarray*}
d_{q_1}(1) & = & 0011011011  \\
d_{q_2}(1) & = & 00111101101011.
\end{eqnarray*}

We start by reading characters from each string into the ``string read''
    expression.

\[
\begin{array} {lll}
\mathrm{String\ 1} & \mathrm{String\ 2} & \mathrm{String\ read}\\ \hline
001101101011 & 0011110110101011 & \mathrm{empty} \\
01101101011 & 011110110101011 & 0 \\
1101101011 & 11110110101011 & 00 \\
101101011 & 1110110101011 & 001 \\
01101011 & 110110101011 & 0011 
\end{array}
\]

At this point we see that the next characters to read from String~1 
    and String~2 are different.
We use a result that is only observed computationally,
    and has no theoretical reason for being true.
This is that the size of every part that depends on the value 
    of $k$ is of the same size, which is known before the computation
    begins.
So an expression $d_{q_k}(1) = v_1 (w_1)^k v_2 (w_2)^k \cdots$ would have all 
    $|w_i|$ constant, and known in advance.
In this case, we are assuming that this size is 2.
So we check if the next two characters of String 2 are the same as the 
    previous two characters of String 1.
(In this case, both of these are ``11''.)
We then truncate the result to give something of the form 
    $(11)^k$ which is valid for both strings.

So we continue.
\[
\begin{array} {lll}
\mathrm{String\ 1} & \mathrm{String\ 2} & \mathrm{String\ read}\\
\hline
01101011 & 110110101011 & 0011  \\
01101011 & 0110101011 & 00(11)^k \\
1101011 & 110101011 & 00(11)^k0 \\
 \vdots  & \vdots  & \vdots   \\
1 & 011 & 00(11)^k0110101 
\end{array}
\]

Again we check if the next two characters of String 2 are equal to the 
    previous two characters in String 1.
We also notice that the two characters ``01'' are in fact repeated more times
    than this, so we get
\[
\begin{array} {lll}
\mathrm{String\ 1} & \mathrm{String\ 2} & \mathrm{String\ read}\\
\hline
1 & 011 & 00(11)^k0110101  \\
1 & 1 & 00(11)^k01101(01)^k  \\
1 & 1 & 00(11)^k011(01)^{k+1}  \\
\mathrm{empty} & \mathrm{empty} & 00(11)^k011(01)^{k+1}1
\end{array}
\]

So we would conjecture that $d_{q_k}(1) = 00(11)^k011(01)^{k+1}1$.

It should be pointed out that this is in no way a proof that this is the 
    general result.  
This has to be done separately in Sections \ref{sec:verify} and \ref{sec:gl}. 

\subsection{Verifying conjecture}
\label{sec:verify}

In this section we show how, given a conjectured expansion for $q_k$, 
    we would verify that this is a valid expansion for all $q_k$.
It should be noticed that this does not prove what type of $\beta$-expansion
    it is (i.e., greedy, lazy ...).
This will be done in Section \ref{sec:gl}.

We will demonstrate this method, by considering an example. 
Consider the greedy expansion $d_{\beta_k}(1) = 1(101)^k1 = 11(011)^k$ 
    associated with 
    the greedy expansion of $1$ for the Pisot root associated with
    \[ P_k^*(x) = P_{\psi_2}(x) x^{3 k +2} + B_{\psi_2}(x)
      = (x^3-x^2-x-1) x^{3 k + 2} + (x+1).\]
For convenience we write $\beta$ for this root (where $\beta$ will depend on $k$).
We see then that this expansions implies 
\[
\frac{1}{\beta} + \frac{1}{\beta^2} + 
    \frac{1}{\beta^4} + \frac{1}{\beta^5} + 
    \frac{1}{\beta^7} + \frac{1}{\beta^8} +  \cdots + 
    \frac{1}{\beta^{3 k + 1}} + \frac{1}{\beta^{3 k+2}} =1 \]
This simplifies to 
$$ 
{\beta}^{-1}+{\beta}^{-2}+
  \left( {\beta}^{-4}+{\beta}^{-5} \right) 
  \sum_{j=0}^{k-1} \left( {\beta}^{3 j} \right)^{-1} = 1.
$$
By subtracting $1$ from both sides, and clearing the denominator,
    this is equivalent to $D_k(\beta) = 0$ where
$$ D_k(x) := - x^{3 k+5} +x^{3 k+4} +x^{3 k + 3} + x^{3 k + 2} - x - 1 $$
But we notice that
    \[D_k(x) = -(P_{\psi_2}(x) x^{3 k +2} + B_{\psi_2}(x)),\]
    hence $D_k(x) = - P_k^*(x)$.
All of these processes can be automated.  
The hardest part is finding a co-factor $C_k(x)$ such that 
    $D_k(x) = C_k(x) P^*_k(x)$.
(We are not always so lucky that $C_k(x) = -1$ as was the case in this example.)
Here we noticed computationally that $C_k(x)$ is always of the 
    form: 
$$
C_k(x) = a_n x^{b_n k + c_n} + a_{n-1} x^{b_{n-1} k + c_{n-1}} +  
                   \cdots + a_2 x^{b_2 k + c_2} + a_{1} x^{b_{1} k + c_{1}}.
$$
For our purposes it was unnecessary to prove that this is always the case,
    as we could easily verify it for all cases checked, and 
    we were using this as a tool to verify the conjectured general form.

\subsection{Check greedy/lazy/univoque/periodic self-bracketed 
    $\beta$-expansion}
\label{sec:gl}

In this section we discuss how one would check if an expression 
    (conjectured using the techniques of Section \ref{sec:conj}
     and verified as a $\beta$-expansion in Section \ref{sec:verify}), 
    is in fact a greedy, lazy or periodic self-bracketing $\beta$-expansion.
Consider a general expression of the form
    \[ E(k) := v_1 (w_1)^{k} v_2 (w_2)^k \cdots (w_{n-1})^k v_n 
       (u_1 (w_{n})^k \cdots (w_{n+m})^k u_m)^\infty \]
    where the $w_i$ all have the same length (this is in fact the case for
    all problems that we studied).
Then the main thing to notice is that there exists a $K$ such that if the 
    $\beta$-expansion $E(K)$ has a desired property (either being or not being
    greedy, lazy, etc), then for all $k \geq K$ we have $E(k)$ has
    the same property.
Moreover the $K$ is explicitly computable, being a function of the 
    lengths of the $v_i$, $w_i$ and $u_i$.
This means that what initially looks like an infinite number of calculations 
    is in fact a finite number of calculations.
The way to see this is that for sufficiently large $k$, most of the 
    comparisons will be done between the $w_i$'s, and then an
    increase in $k$ will not change this, but just add another redundant 
    check to something already known.

\bigskip
\bigskip

\noindent
\section{Comments, Open Questions and Further Work}

There are some interesting observations that can be made from the data
    and results so far.  
This investigation has opened up a number of questions.
\begin{itemize}
\item First,
    given a sequence of greedy of lazy $\beta$-expansions
    of a nice sequence of Pisot numbers $q_k$ that looks like:
    \[ E(k) := v_1 (w_1)^{k} v_2 (w_2)^k \cdots (w_{n-1})^k v_n 
       (u_1 (w_{n})^k \cdots (w_{n+m})^k u_m)^\infty \]
    is it always true that $|w_1| = |w_2| = \cdots = |w_{n+m}|$?
\item 
    Is the co-factor from Section \ref{sec:verify} always of the form:
$$
C_k(x) = a_n x^{b_n k + c_n} + a_{n-1} x^{b_{n-1} k + c_{n-1}} +  
                   \cdots + a_2 x^{b_2 k + c_2} + a_{1} x^{b_{1} k + c_{1}}?
$$
\item It appears in Table \ref{tab:SalemUnivoque} that for all Salem numbers
of degree $4$ and $6$
    greater than $\approx 1.83$, these Salem numbers are univoque.
   Is this just an artifact of small degrees, or is something more 
    general going on?
\item In general, are the greedy/lazy $\beta$-expansions even periodic 
    for Salem numbers?
    (This is not known to be true, see \cite{Boyd96a} for more details.)
\item It is known that Pisot numbers can be written as a limit of Salem numbers, 
    where if $P(x)$ is the minimal polynomial of a Pisot number, then 
    $P(x) x^n \pm P^*(x)$ has a Salem number as a root, which tends to the root of the
    Pisot number.  
    Some preliminary and somewhat haphazard investigation suggests that we 
    might be able to find a ``regular'' looking expression for the greedy 
    (resp.~lazy) $\beta$-expansion of these
    Salem numbers, which tends towards the greedy (resp.~lazy) $\beta$-expansion 
    of the Pisot number.  
    If true, this could have implications towards questions concerning 
    the $\beta$-expansions of Salem numbers being eventually periodic.
\end{itemize}

\section*{Acknowledgments}

The authors wish to thank David Boyd for stimulating discussions, and the
referee for a careful reading of the manuscript.

\section*{Note added on June 12, 2006}

Just before submitting this paper we came across a paper where the
topological structure of the set ${\U}$ and of its (topological)
closure are studied. We cite it here for completeness:

\noindent
{\small V. Komornik, P. Loreti, On the structure of univoque sets,
{\em J. Number Theory}, to appear.}

\noindent
One can also read consequences of the results of that paper in

\noindent
{\small M. de Vries, Random $\beta$-expansions, unique expansions and Lochs' 
Theorem, PhD Thesis, Vrije Universiteit Amsterdam, 2005.} 

\noindent
(available at {\tt http://www.cs.vu.nl/$\sim$mdvries/proefschrift.pdf}).


\end{document}